\documentclass{article}
\usepackage{spconf,amsmath,graphicx}

\usepackage{epsfig,psfrag}
\usepackage{pst-all}
\usepackage{amsmath,amsthm,amssymb,amsfonts,upref,cite,epsf,color,bm}
\usepackage{graphicx}
\usepackage{color}
\usepackage{amsmath,mathtools}
\usepackage{graphicx}
\usepackage{calc}

\newtheorem{theorem}{Theorem}[section]

\usepackage{etoolbox}
\makeatletter
\appto{\equation}{
  \notag
  \preto{\label}{
    \refstepcounter{equation}
    \tag{\arabic{equation}}
  }
}
\makeatother

\usepackage{framed}

\newcommand\defeq{:=}

\usepackage[applemac]{inputenc}

\def \expect {{\rm E} }

\def \twiddle[#1] {e^{-j \frac{2 \pi}{N}  #1 }}
\def \twiddleneg[#1] {e^{j \frac{2 \pi}{N}  #1 }}
\DeclareMathOperator*{\argmax}{argmax}

\def\ML_est{\hat{\mathbf{x}}_{\text{ML}}}

\newcommand{\autocovfunc}{\mathbf{R}} 
\newcommand{\ACF}{\autocovfunc} 
\newcommand{\cig}{\mathcal{G}}
\newcommand{\cigrule}{\widehat\cig}

\newcommand{\mincoh}{\rho_{\min}}
\newcommand{\SDM}{\mathbf{S}}

\newcommand{\be}{\begin{equation}}
\newcommand{\ee}{\end{equation}}

\newcommand{\lagvar}{m}

\newcommand{\timeidx}{n} 
\newcommand{\freqidx}{f} 
 
\newcommand{\ensembleidx}{i} 
\newcommand{\samplesize}{N}

\newcommand{\coefflen}{p}

\newcommand{\xobs}{\mathbf{X}}%X_{1}^{\samplesize}}
\newcommand{\xEobsVec}{\tilde{\mathbf{x}}}
\newcommand{\dmax}{d_{\max}}
\newcommand{\processclass}{\mathcal{M}(p,\dmax,\mincoh)}

\allowdisplaybreaks

\makeatother

\usepackage{setspace}
\title{On the Information-theoretic Limits of Graphical Model Selection for
Gaussian Time Series}
\name{Gabor Hannak, Alexander Jung, Norbert Goertz\thanks{Funded by the EU
FP7 project NEWCOM\# (GA 318306).}} \address{Institute of Telecommunications\\
Vienna University of Technology, Austria \\
Email: \{ghannak,  ajung, ngoertz\}@nt.tuwien.ac.at}

\begin{document}
% \tracingall \ninept
\maketitle
\begin{abstract}
We consider the problem of inferring the conditional independence
graph (CIG) of a multivariate stationary dicrete-time Gaussian random process
based on a finite length observation.
Using information-theoretic methods, we derive a lower bound on the error
probability of any learning scheme for the underlying process CIG. This bound,
in turn, yields a minimum required sample-size which is necessary for any
algorithm regardless of its computational complexity, to reliably select the
true underlying CIG. Furthermore, by analysis of a simple selection scheme, we
show that the information-theoretic limits can be achieved for a subclass of
processes having sparse CIG. We do not assume a parametric model for the
observed process, but require it to have a sufficiently smooth spectral density
matrix (SDM).

\end{abstract}

\vspace{-3mm} 
\vspace*{-2mm} 
\section{Introduction}\label{sec:intro}
\vspace*{-2mm}

We consider multivariate time series, i.e., vector-valued discrete
time stationary random processes
$\mathbf{x}[\timeidx]\!\!=\!\!(x_{1}[\timeidx],\ldots,x_{p}[\timeidx])^T \,,
n\in\mathbb{Z}$.
An important problem within multivariate time series analysis, e.g. in the context
of medical data or environmental monitoring data, is the characterization of the
interdependencies between the individual process components
\cite{Brockwell91,Luetkepol2005}. A particular representation of the statistical
relations governing the process components is obtained via the conditional
independence graph (CIG) \cite{Dahlhaus2000}.
In particular, the CIG associated with the time series $\mathbf{x}[\timeidx]$ is
an undirected graph with node set $V=\left\{1,\ldots,p\right\}$, where a
specific node $r$ represents the component process $x_r[\timeidx]$ and an edge
$(r,r')$ represents a dependency between process components $x_r[n]$ and
$x_{r'}[n]$.

By graphical model selection, we refer to the problem of determining the
underlying CIG based on a finite-length observation of the time series.
 The problem of graphical model selection for the case of Gaussian Markov random
 fields (GMRF), considered in \cite{WangWain2010}, is contained in our setting
 as the special case of a white (i.i.d.) process. Thus, our work can be regarded
 as a
generalization of \cite{WangWain2010} in the sense that we allow for temporal
correlation of the observed samples.\\
Numerous graphical model selection schemes for
vector-valued stationary processes have been proposed
\cite{PetersCausalInf2012,Bento2010,Songsiri09,songsiri2010,JuHeck2014}.
Most of them are based on finite-dimensional parametric models. In particular,
\cite{Bento2010,Songsiri09,songsiri2010} model the observed process as a
vector-valued autoregressive (VAR) process.
Recently, the authors of \cite{JuHeck2014} proposed a nonparametric selection
scheme which does not rely on a parametric process model but requires only certain
smoothness constraints to be satisfied by the process spectral density matrix
(SDM). A performance analysis for the proposed scheme in
\cite{JuHeck2014} provides sufficient conditions on the observed
sample-size that enables reliable model selection. 

In this work we complement
the sufficient conditions presented in \cite{JuHeck2014} with
lower-bounds on the sample-size required by any graphical model selection scheme
to be reliable. We highlight that, similar to \cite{WangWain2010}, these
necessary conditions apply to any model selection scheme regardless of its
computational complexity.

Our results apply to the high-dimensional regime,
where the process dimension $p$ and the sample-size $\samplesize$ may diverge
simultaneously. In particular, we allow for scenarios where the sample-size $N$
is much smaller than the process dimension $p$. To allow for accurate selection
schemes in the high-dimensional regime, one needs to require additional structural properties of the process.
The structure induced by requiring the CIG to be sparse allows for reliable
selection even in the high-dimensional regime. Our results can be used as a
standard against which the performance (in terms of required sample-size) of a
specific selection scheme can be compared.\\
We begin in Section \ref{sec:problem_formulation} with a discussion of the
set-up and the problem formulation. In Section \ref{sec:main_results} we
introduce the main results of this work, Theorem \ref{thm_gmrf} and Theorem
\ref{thm_processes}, and discuss their consequences. In Section \ref{proofs} we
give a high-level outline of the proofs for the main results.
% In Section \ref{sec:conclusions} we conclude.

\textbf{Notation.}  The identity matrix of dimension $L\times L$ is denoted
$\mathbf{I}_{L}$, without the subscript if the dimension is clear by
context. We denote the $k$th column of the identity matrix $\mathbf{I}$ by
$\mathbf{e}_{k}$. Given a matrix $\mathbf{A}$, we denote its entry in the $i$th
row and $j$th column by $\left(\mathbf{A}\right)_{ij}$. The set $\left\{
1,\ldots,p\right\}$ is denoted by $[p]$. The Kronecker product of the matrices
$\mathbf{A}$ and $\mathbf{B}$ is denoted $\mathbf{A} \otimes \mathbf{B}$.
We define the Kronecker delta $\delta[\timeidx]$ as $\delta[0] = 1$ and
$\delta[\timeidx]=0$ otherwise. Given an index set $\mathcal{A}$, we denote
$\mathbf{I}_{\mathcal{A}} = \sum_{r \in \mathcal{A}} \mathbf{e}_{r}
\mathbf{e}_{r}^{T}$ and $\mathbf{1}_{\mathcal{A}} = \sum_{r \in \mathcal{A}}
\mathbf{e}_{r}$, respectively.

%%%%%%%%%%%%%%%%%%%%%%%%%%%%%%%%%%%%%%%%%%%%%%%%%%%%%%%%%%%%%%%%%%%%%%%%%%%%%
%%%%%%%%%%%%%%%%%%%%%%%%%%%%%%%%%%%%%%%%%%%%%%%%%%%%%%%%%%%%%%%%%%%%%%%%%%%%%

\vspace*{-3mm}
\section{Problem Formulation}\label{sec:problem_formulation}
\vspace*{-2mm}
% \vspace*{-2mm}

Consider a zero-mean multivariate Gaussian time series
$\mathbf{x}[\timeidx]\!\!=\!\!\left(x_1[\timeidx],\ldots,x_p[\timeidx]
\right)^T$, $n\!\in\!\mathbb{Z}$.
The temporal dependence between the samples $\mathbf{x}[\timeidx]$ is captured
by the matrix-valued autocorrelation function (ACF)
\vspace{-2mm}
\begin{equation}
\autocovfunc_x\!\left[\lagvar\right]\!
=\mathrm{E}\left\{\mathbf{x}\!\left[\lagvar\right]\mathbf{x}^T\!\left[0\right]\right\}\vspace{-2mm}
\,,
\end{equation}
assumed to be summable, i.e., $\sum_{k=-\infty}^{\infty}
\|\ACF_x\!\left[\lagvar\right]\|_{\infty}\!<\!\infty$. The spectral
density matrix (SDM),
\vspace{-2mm}
\begin{equation}\label{eq:sdm}
\SDM_x\!\left(\theta\right)\!
=\sum_{m=-\infty}^{\infty} \!\! \ACF_x\!\left[\lagvar\right] \exp(-j2\pi\theta
\lagvar) \,,\vspace{-2mm}
\end{equation}
describes the correlation structure between the process components in the
frequency domain. In what follows, we assume
\vspace{-2mm}
\begin{equation}\label{eq:AB}
1 \leq \lambda_{\min}(\SDM_x(\theta)) \leq
\lambda_{\max}(\SDM_x(\theta)) \leq \mathrm{B} < \infty \,,\vspace{-1mm}
\end{equation}
where $\lambda_{\min}(\SDM_{x}(\theta))$ and $\lambda_{\max}(\SDM_{x}(\theta))$
denote the minimum and maximum eigenvalue of  the SDM, respectively. The lower bound in \eqref{eq:AB} ensures certain Markov properties of the CIG
\cite{LauritzenGM} and the upper bound follows from the summability of the ACF.

Our analysis applies to processes with a smooth SDM, i.e., the entries of
$\SDM_{x}(\theta)$ are smooth functions. Due to the Fourier relationship
(\ref{eq:sdm}), these smoothness constraints can be expressed via the ACF moment
\vspace{-2mm}
\begin{equation}\label{equ_def_moment_x}
\mu_x = \sum_{\lagvar=-\infty}^{\infty} \lvert \lagvar \rvert
\|\ACF_x\!\left[\lagvar\right]\|_{\infty} \,.\vspace{-2mm}
\end{equation}
For a small moment $\mu_x$, the ACF has to be well concentrated around
$\lagvar\!=\!0$.

Given a process $\mathbf{x}[\timeidx]$, we define its underlying conditional
independence graph (CIG) as $\cig\!=\!(V,E)$ with node set $V\!=\![p]$ and edge
set $E\!\subseteq\! V\! \times\! V$.
The nodes of $\mathcal{G}$ represent the scalar component processes
$x_r[\timeidx]$.
The edge set is characterized by requiring $(r,r')\notin E$ if and only if the
process components $x_r[\timeidx]$ and $x_{r'}[\timeidx]$ are conditionally
independent given the remaining components $\left\{x_t[\timeidx] \mid t\!\in\!
\left[p\right]\!\setminus\!\left\{ r,r' \right\}\right\}$. For a Gaussian
process with SDM satisfying (\ref{eq:AB}), the CIG can be characterized
conveniently in terms of the SDM. In particular
\cite{Brillinger96remarksconcerning},
\vspace{-2mm}
\begin{equation}
\label{equ_charc_edgeset_Gaussian_process}
(r,r')\!\notin\! E \mbox{, if and only if } (\SDM_x^{-1}(\theta))_{rr'}\!=\!0
\ \quad \forall \theta\!\in\![ 0,1).
\end{equation}
To quantify the strength of the dependence
between connected process components, we define the minimum partial spectral coherence of 
a process $\mathbf{x}[\timeidx]$ as
\vspace{-2mm}
\begin{equation}
\label{equ_def_rho_x}
\rho_{x} = \min_{(r,r') \in E} \left( \int_{0}^{1} 
\frac{\left(\SDM_x^{-1}(\theta)\right)_{rr'}^2}{\left(\SDM_x^{-1}(\theta)\right)_{rr}\left(\SDM_x^{-1}(\theta)\right)_{r'r'}}
d\theta \right)^{\frac{1}{2}} \,.\vspace{-2mm}
\end{equation}
The neighborhood and degree of a node $r$ are defined as
$\mathcal{N}(r)\!=\!\left\{ r' \mid (r,r')\!\in\! E\right\}$, and as
$d_r\!=\!\lvert \mathcal{N}(r) \rvert$, respectively.
We consider processes whose CIG have bounded degrees, i.e., for some (typically
small) $d_{\max}$
\vspace{-2mm}
\begin{equation}\label{defdmax}
d_{r} \leq d_{\max} \quad \forall r\in [p] \,. \vspace{-2mm}
\end{equation}
In the following $\processclass$ denotes the class of $p$
dimensional vector processes whose underlying CIGs have maximum
degree at most $\dmax$, whose SDMs fulfill (\ref{eq:AB}) with
some constant $\mathrm{B}$ and have minimum partial coherence not smaller than
$\mincoh$.\\
A graphical model selection scheme $\cigrule$ maps the observed samples
$\xobs\!=\!\left(\mathbf{x}[1],\ldots,\mathbf{x}[\samplesize]\right)
\in \mathbb{R}^{p \times \samplesize}$ to an estimate $\cigrule$ of the true
CIG $\cig$.
We define the maximum selection error probability of a scheme
$\cigrule(\cdot)$ as
\vspace{-2mm}
\begin{equation}\label{def:error_probability}
p_{err}(\cigrule) = \max_{x\in \mathcal{M}} \Pr \left\{
\cigrule(\xobs) \neq \cig \right\} \,. \vspace{-2mm}
\end{equation}
For a specific sample-size N, we define
the minimax detection error as
\vspace{-2mm}
\begin{equation}
p_{err}(\samplesize) = \min_{\cigrule} p_{err}(\cigrule) \,.\vspace{-2mm}
\end{equation}
We then study conditions on the sample-size
$\samplesize(p,\dmax,\mincoh)$ as a function of the remaining problem
parameters, such that asymptotically reliable selection is achievable, that is,
$p_{err}(\samplesize)\!\rightarrow \!0$ as
$\samplesize(p,\dmax,\mincoh)\! \rightarrow\! \infty$.

%%%%%%%%%%%%%%%%%%%%%%%%%%%%%%%%%%%%%%%%%%%%%%%%%%%%%%%%%%%%%%%%%%%%%%%%%%%%%
%%%%%%%%%%%%%%%%%%%%%%%%%%%%%%%%%%%%%%%%%%%%%%%%%%%%%%%%%%%%%%%%%%%%%%%%%%%%%
\vspace{-2mm}
\section{Main Results}\label{sec:main_results}
\vspace{-2mm}

\subsection{Necessary Conditions for Consistent Model Selection}
For our first result, we closely follow the argument in \cite{WangWain2010}
used to derive lower-bounds on the required sample-size for i.i.d.~samples.
The observation of i.i.d.~samples is contained in our setup as the special case
obtained for $\mu_x \!=\! 0$, implying that the SDM $\SDM_x(\theta)$ is flat,
i.e., it does not depend on $\theta$.

\begin{theorem}\label{thm_gmrf}
Consider a process in the class $\processclass$
with $\mincoh\!\in\!(0,\frac{1}{4}]$. A necessary condition for
asymptotically reliable graphical model selection is
\begin{equation}\label{necessary_cond_thm_iid}
\samplesize>\frac{\log{p \choose 2}-1}{4\mincoh^2} \,.
\end{equation}
\end{theorem} 
\noindent
A proof is sketched in Section \ref{proofs}.

The lower bound \eqref{necessary_cond_thm_iid} on the sample-size completely
ignores the correlation width of the process, which is quantified by
$\mu_{x}$.
Intuitively, we would expect that with increasing correlation width the required
sample-size becomes larger. This is also reflected by the sufficient conditions
on the sample-size in \cite{JuHeck2014} for a novel nonparametric selection scheme. In
particular the results of \cite{JuHeck2014} suggest that the sample-size has to grow
proportionally to $\mu^{2}_{x}$.
% 
% The case of i.i.d.~samples represents an optimum situation in the sense that the
% information gain by observing an additional sample is maximal.
% Consider a process $\mathbf{x}[\timeidx]$ with ACF $\ACF_x[\lagvar]=0$ for
% $ \lvert m\rvert \geq K$ and a selection algorithm which requires $L$
% i.i.d.~samples.
% Then, if $KL$ samples of $\mathbf{x}[\timeidx]$ are observed, one can
% extract $L$ i.i.d~samples by simple subsampling the observed
% block by a factor $K$.

Another argument supporting the intuition that that the required sample-size $N$
has to be larger for increasing $\mu_x$ is that for larger values of $\mu_{x}$
the SDM has faster variations over $\theta$. Therefore, in order to determine
the zeros of the inverse SDM, we have to estimate the SDM values at more
sampling points (placed denser). However, the values of the
SDM at different frequencies $\theta$ are strongly coupled via the condition
\eqref{equ_charc_edgeset_Gaussian_process}. This implies via \eqref{defdmax}
(with a small $d_{\max}$) that the values of the inverse SDM $\SDM^{-1}(\theta)$
have a small joint support, which corresponds to the edge set E of the CIG.
\vspace{-4mm}

\subsection{Sufficient Condition for Consistent Model Selection}

We next show that, at least for the special case of $d_{\text{max}}=1$, the
coupling via \eqref{equ_charc_edgeset_Gaussian_process} compensates the effect
of increasing correlation width $\mu_{x}$ such that the required sample-size is
independent of $\mu_{x}$. This will be accomplished by analyzing a specific
model selection scheme. This scheme is similar to the exhaustive search decoder
used in \cite{Wain2009TIT} for the derivation of sufficient conditions for
sparsity recovery in the high-dimensional sparse linear model.

In order to keep the argument as simple as possible, we assume that
\begin{itemize}
\item the ACF is exactly supported within
$\{-\frac{\samplesize}{2}+1,\ldots,0,\ldots,\frac{\samplesize}{2}-1\}$, i.e.,
\vspace{-2mm}
\begin{equation}
\label{equ_suff_cond_ass_finite_support_ACF}
\autocovfunc_{x}[\lagvar] = 0 \mbox{ if }  |\lagvar| \geq \frac{\samplesize}{2}
\,,\vspace{-2mm}
\end{equation} 
\item the ACF is real-valued and symmetric (instead of being merely Hermitian
symmetric), i.e.,
\vspace{-2mm}
\begin{equation} 
\label{equ_suff_cond_ass_symmetric_ACF}
\autocovfunc_x[-\lagvar] = \autocovfunc_{x}[\lagvar] \,. \vspace{-2mm}
\end{equation} 
\end{itemize}
A specific subclass of processes satisfying these assumptions is obtained by
applying a real-valued scalar filter component-wise to a white noise vector
process, i.e.,
\vspace{-2mm}
\begin{equation}\label{proc_filter}
\mathbf{x}[\timeidx] = \sum_{\lagvar=0}^{K-1} h[\lagvar] \mathbf{w}[\timeidx -
\lagvar],\vspace{-2mm}
\end{equation}
with $h[\lagvar]$ being a length-$K$ filter impulse response with $K <
\samplesize/4$. We assume that the filter $h[\timeidx]$ is normalized such that
$\sum_{m=-\infty}^{\infty} h^{2}[\timeidx] = 1$, which by Parseval's
theorem implies $\sum_{\freqidx=0}^{2\samplesize-1}
\lvert H(\frac{\freqidx}{2\samplesize})\rvert^{2} = 2\samplesize$, where
$H(\theta) = \sum_{\timeidx=-\infty}^{\infty} h[\timeidx] \exp(-j2\pi \timeidx
\theta)$ denotes the discrete time Fourier transform of the impuse response $h[\timeidx]$.

Here, $\mathbf{w}[\timeidx]$ is an i.i.d. zero-mean Gaussian process with
marginal covariance matrix $\mathbf{C}$, i.e., $\mathbf{w}[\timeidx] \sim
\mathcal{N}(\mathbf{0},\mathbf{C})$. The associated precision matrix
$\mathbf{C}^{-1}$ has at most one non-zero off-diagonal entry in each row. This
implies that the CIG of $\mathbf{x}[\timeidx]$ satisfies \eqref{defdmax} with
$d_{\max}=1$.

The constraint \eqref{equ_suff_cond_ass_finite_support_ACF} is satisfied for
processes with a smooth SDM or, equivalently, a small moment $\mu_{x}$.
Typically the ACF will not be exactly zero for $|m| \geq \samplesize/2$. We make this
idealized assumption only to keep our argument as simple as possible.  However, we expect
that our main conclusions are also valid for any process with sufficiently small
correlation width, i.e., small $\mu_{x}$.

Following a method in \cite{CaiZhaoZhou2013}, we define, given the observations
$\mathbf{x}[1],\ldots,\mathbf{x}[\samplesize+1]$, the enlarged observation set
$\xEobsVec[1],\ldots,\xEobsVec[2 \samplesize]$ by
\vspace{-2mm}
\begin{equation}
\xEobsVec[\timeidx] = \begin{cases} \mathbf{x}[\timeidx]  &\mbox{for
} \timeidx \in [\samplesize+1] \\  \hspace*{-1mm} \mathbf{x}[2\samplesize
\!-\!\timeidx\!+\!2]  &\mbox{for } \timeidx \in [2\samplesize] \setminus
[\samplesize+1]. \end{cases}\vspace{-2mm}
\end{equation}
Because of \eqref{equ_suff_cond_ass_finite_support_ACF} and
\eqref{equ_suff_cond_ass_symmetric_ACF}, the covariance matrix
$\widetilde{\mathbf{C}}$ of the stacked vector $\xEobsVec \defeq \left(
\xEobsVec^{T}[1],\ldots,
\xEobsVec^{T}[2\samplesize] \right)^{T}$ can be shown to be a block
circulant matrix \cite{InverseBCM83} with first row given by
\vspace{-2mm}
\begin{equation}
\label{equ_def_first_row_cov_matrix_tilde_x}
\left( \autocovfunc_x[0], \ldots, \autocovfunc_x[\samplesize],  \autocovfunc_x[\samplesize-1],\ldots, \autocovfunc_x[2] \right).\vspace{-2mm}
\end{equation}

Let us define the DFT of the enlarged observation set as
\vspace{-2mm}
\begin{equation}
\hat{\mathbf{x}}[\freqidx] \defeq (1/\sqrt{2\samplesize}) \sum_{\timeidx \in [2
\samplesize]} \xEobsVec[\timeidx] \exp\left( -j 2 \pi
\frac{(\freqidx-1) (\timeidx-1)}{2\samplesize}  \right) \vspace{-2mm}
\end{equation}
Some calculation reveals that the covariance matrix $\widehat{\mathbf{C}}$ of
the stacked vector $\hat{\mathbf{x}} \defeq \left(
\hat{\mathbf{x}}^{T}[1],\ldots, \hat{\mathbf{x}}^{T}[2\samplesize] \right)^{T}$
satisfies
\vspace{-2mm}
\begin{equation}
\widehat{\mathbf{C}} =
\frac{1}{2\samplesize}\mathbf{P}\widetilde{\mathbf{C}}
\mathbf{P}^{H},\vspace{-2mm}
\end{equation} 
with the matrix $\mathbf{P}$ as defined in \cite[p. 809]{InverseBCM83}. Since
$\widetilde{\mathbf{C}}$ is block circulant, it follows that
$\widehat{\mathbf{C}}$ is a block diagonal matrix with the $\freqidx$th
$\coefflen \times \coefflen$ diagonal block given by the $\freqidx$th bin of the
DFT of the row \eqref{equ_def_first_row_cov_matrix_tilde_x}, which can be shown
to coincide with $\SDM_{x}((\freqidx-1)/2\samplesize)$.
Therefore, the vectors $\hat{\mathbf{x}}[\freqidx]$ are independent (across
$\freqidx$) zero-mean Gaussian vectors with covariance matrix
$\SDM_{x}((\freqidx-1)/2\samplesize)$. 

Let $\hat{\mathbf{x}}_r$ denote the $r$th
row of $\left( \hat{\mathbf{x}}[1],\ldots, \hat{\mathbf{x}}[2\samplesize]
\right)$. The estimation of the neighborhood $\mathcal{N}(r)$, under the
assumption that the true CIG satisfies \eqref{defdmax} with $d_{max}=1$, is
carried out for every $r\in V$ by:
\begin{itemize}
\item for all $r' \in [\coefflen] \setminus \{r\}$ compute the statistic
\vspace{-2mm}
\begin{equation}\label{decoder1}
Z(r') \defeq \hat{\mathbf{x}}_r^T \hat{\mathbf{x}}_{r'} \,,\vspace{-2mm}
\end{equation}
\item determine the maximizing index
\vspace{-2mm}
\begin{equation}\label{decoder2}
\hat{r} \defeq \argmax_{r' \in [\coefflen] \setminus \{r\}} |Z(r')| \,,\vspace{-2mm}
\end{equation}
\item compare the maximum statistic with threshold $\eta$ to obtain  
\begin{equation}\label{decoder3}
\widehat{\mathcal{N}}(r) = \begin{cases} \left\{\hat{r}\right\} & \mbox{if }
|Z(\hat{r})| \geq \eta \\ \emptyset & \mbox{otherwise.} \end{cases}
\end{equation}
\end{itemize}
The choice of $\eta$ will be discussed in Section \ref{proof_processes}.
Note that this algorithm does not necessarily produce a CIG satisfying
\eqref{defdmax} with $d_{\max}=1$.

The following result shows that for the subclass of
$\mathcal{M}(p,\dmax=1,\mincoh)$ given by processes of the form
\eqref{proc_filter}, the simple selection scheme
\eqref{decoder1}-\eqref{decoder3} achieves the information theoretic limit
stated in Theorem \ref{thm_gmrf}.

\begin{theorem}\label{thm_processes}
Consider a process in the class $\mathcal{M}(p,d_{\max}=1,\mincoh)$ generated
according to \eqref{proc_filter}. A sufficient condition on the sample-size $N$
for achieving probability of incorrect selection not larger than $\delta$ using
the selection scheme given by \eqref{decoder1}-\eqref{decoder3} is
\vspace{-2mm}
\begin{equation}\label{bound_processes}
\samplesize>
\frac{32B^4}{\mincoh^2} \log \left(\frac{2p^2}{\delta} \right) \,.
\end{equation}
\end{theorem}
\noindent
A proof is sketched in Section \ref{proofs}.\\
\textbf{Discussion.} The bound \eqref{bound_processes} matches the scaling of
the necessary condition \eqref{necessary_cond_thm_iid} in Theorem
\ref{thm_gmrf}, in particular, it scales proportional to $\frac{1}{\mincoh^2}$
and to $\log p$. We deduce that for the processes of the form \eqref{proc_filter} in the subclass
$\mathcal{M}(p,d_{\max}=1,\mincoh)$ the theoretic limit stated in Theorem
\ref{thm_gmrf} can be achieved. Thus, since \eqref{bound_processes} does not
depend on $\mu_{x}$, for this special case the temporal correlation of the
process, quantified by $\mu_{x}$, does not increase the required sample-size for reliable graphical
model selection.

\vspace*{-2mm}
\section{Proof Sketches}\label{proofs}
\vspace*{-2mm}

\subsection{Proof of Theorem \ref{thm_gmrf}}
\vspace{-2mm}
For the derivation of Theorem \ref{thm_gmrf}, we closely follow the method put
forward in \cite{WangWain2010}.
In particular, the result is based on a finite ensemble $\mathcal{M}_0$
containing $M\in\mathbb{N}$ different processes
$\mathbf{x}^{(i)}[\timeidx]\!\in\!\processclass$ with associated SDM
$\mathbf{S}_{i}(\theta)$ and CIG $\cig(i)$, respectively.

Assuming that the observed process $\mathbf{x}[\timeidx]$ is taken uniformly at
random out of $\mathcal{M}_{0} \subseteq \processclass$, we may interpret the
graphical model selection problem as a
communication problem:
using a random index $i$, distributed uniformly over the set $[M]$, we  
select the process $\mathbf{x}^{(i)}[\timeidx]$ as the observed process,
i.e., $\mathbf{x}[\timeidx] \!=\! \mathbf{x}^{(i)}[\timeidx]$. Based on the
samples $\xobs=(\mathbf{x}[1],\ldots,\mathbf{x}[N])$, the problem of selecting
the true graphical model is now equivalent to detecting $\cig(i)$.
A selection rule $\cigrule(\xobs)$ can be interpreted as a decoder, mapping the
observation $\xobs$ to an estimate of the true CIG $\cig(i)$.

The maximum probability $p_{err}(\hat{G})$ (cf. \eqref{def:error_probability})
probability \eqref{def:error_probability} of any selection rule can be lower bounded via Fano's inequality \cite{coverthomas} as
\vspace{-2mm}
\begin{equation}\label{eq:fano}
p_{err}(\widehat{\mathcal{G}}) \geq 1 - \frac{\mathrm{I}(\xobs;i)+1}{H(\cig(i))} \,.\vspace{-2mm}
\end{equation}
Here, $H(\cig(i))$ denotes the entropy of the random CIG $\cig(i)$ associated
with the process $\mathbf{x}^{(i)}[\timeidx]\!\in\!\mathcal{M}_{0}$, which is
selected uniformly at random from $\mathcal{M}_{0}$.
Since the bound in \eqref{eq:fano} applies to any selection rule, it is also a
lower bound on the minimax error probability $p_{err}(\samplesize)$, i.e.,
$p_{err}(\samplesize) \geq 1 - \frac{\mathrm{I}(\xobs;i)+1}{H(\cig(i))}$.
Therefore, for asymptotically reliable model selection, i.e.,
$\lim_{N\rightarrow\infty} p_{err}(\samplesize) = 0 $, a necessary condition is
\vspace{-2mm}
\begin{equation} 
\label{equ_cond_MI_asmpt_reliable_sel}
I(\xobs;i) \geq H(\cig(i)) - 1.\vspace{-2mm}
\end{equation} 
Based on \eqref{equ_cond_MI_asmpt_reliable_sel}, we will now derive necessary
conditions on the sample-size $\samplesize$ by using upper bounds on the mutual
information $I(\xobs;i)$ which depend explicitly on $\samplesize$. In
particular, since given the index $i$, the observation $\mathrm{vec}(\xobs)$ is a multivariate normal
vector with zero-mean and covariance matrix $\mathbf{C}_{i} \defeq \expect \{
\mathrm{vec}(\xobs) \mathrm{vec}(\xobs)^{T} | i \}$, we can use the following
\emph{entropy-based} upper bound \cite{WangWain2010,WeiWangPhd}
\vspace{-2mm}
\begin{equation}
\label{equ_entropy_based_bound_MI}
I(\xobs;i) \leq \log \big| \overline{\mathbf{C}} \big| - (1/M) \sum_{i\in [M]}
\log \big|\mathbf{C}_{i} \big|,\vspace{-2mm}
\end{equation} 
where $\overline{\mathbf{C}} \defeq (1/M) \sum_{i\in [M]}  \mathbf{C}_{i}$.

%We derive the necessary conditions on the sample-size one obtains from the
%right side of (\ref{eq:fano}) going to $0$. To this end we will make use of the
%following inequality
%\begin{align}
%I(i,\mathbf{X}_1^N) =& h(\mathbf{X}_1^N) - h(\mathbf{X}_1^N\vert i) \\
%\leq & \log\det\cov(\mathbf{X}_1^N) \\
%& - \frac{1}{M} \sum_{m=1}^{M} \log \det \cov (\mathbf{X}_1^N \vert i=m) \,.
%\end{align}

% The set of integers $\left\{1,\ldots,p \right\}$ is denoted by $[p]$.
 Given $p$ nodes, let $\mathcal{S}(\ensembleidx)$ denote an
enumeration of the $\bar{p} = \binom{p}{2}$ different simple graphs (with node
set $[p]$) containing a single edge. Let us define a bijective map
$\mathcal{S}(i)$ which assigns an edge $(r_{\ensembleidx},r'_{\ensembleidx})$
uniquely to an index $\ensembleidx \in [\bar{p}]$. Consider the ensemble
$\mathcal{M}_{0}$ of size $M=\bar{p}$, constituted by the Gaussian processes
$\mathbf{x}^{(i)}[\timeidx]$ with SDM
\vspace{-2mm}
\begin{equation} 
\label{equ_def_SDM_ensemble_i_i_d}
\mathbf{S}_{i}(\theta) =   2\mathbf{I} -
\frac{2\mincoh}{1+4\mincoh}\mathbf{1}_{\mathcal{S}(i)}\mathbf{1}_{\mathcal{S}(i)}^T \vspace{-2mm}
\end{equation}
for $i\in[\bar{p}]$.
Note that the SDM in \eqref{equ_def_SDM_ensemble_i_i_d} does not depend on
$\theta$. Therefore, due to the Fourier relationship \eqref{eq:sdm}, the
corresponding ACF is given by $\mathbf{R}_{i}[\lagvar] = \mathbf{S}_{i}(0) 
\delta[\lagvar]$, implying that $\mu_{x} = 0$ for all $\mathbf{x}[\timeidx] \in
\mathcal{M}_{0}$ (cf. \eqref{equ_def_moment_x}).\\
The CIG $\cig(i)$ associated with \eqref{equ_def_SDM_ensemble_i_i_d} contains a
single edge between nodes $r_{i}$ and $r_{i'}$. Therefore,
\eqref{defdmax} is satisfied for any $\dmax \geq 1$.
Moreover, each process $\mathbf{x}^{(i)}[\timeidx] \in \mathcal{M}_{0}$ with SDM
given by \eqref{equ_def_SDM_ensemble_i_i_d} has a unique CIG $\mathcal{G}(i)$,
i.e., $\cig(i) \neq \cig(i')$ for $i \neq i'$ and, in turn,
\vspace{-2mm}
\begin{equation} 
\label{equ_entropy_graph_ensemble_i_i_d}
H(\mathcal{G}(i)) = \log_{2} | M | = \log_{2} \bar{p}. \vspace{-2mm}
\end{equation} 
If $\mincoh \leq 1/4$, the eigenvalues of the SDM in
\eqref{equ_def_SDM_ensemble_i_i_d} satisfy \eqref{eq:AB} (with $B=3$).
Applying the matrix inversion lemma \cite{HendersonSearle} to
\eqref{equ_def_SDM_ensemble_i_i_d}, we obtain
\vspace{-2mm}
\begin{equation} 
\label{equ_exp_inv_SDM_ensemble_i_i_d}
\mathbf{S}_{i}^{-1}(\theta) =  \mathbf{I} + 2\mincoh
\mathbf{1}_{\mathcal{S}(i)}\mathbf{1}_{\mathcal{S}(i)}^T. \vspace{-2mm}
\end{equation}
Using \eqref{equ_exp_inv_SDM_ensemble_i_i_d} and $\mincoh \leq 1/4$, the
minimum partial coherence $\rho_{x}$ (cf. \eqref{equ_def_rho_x}) of the process
$\mathbf{x}^{(i)}[\timeidx]$ can be shown to satisfy
%\begin{align}
%\rho_{x} & =  \left( \int_{0}^{1} 
%\frac{\left(\SDM_x^{-1}(\theta)\right)_{rr'}^2}{\left(\SDM_x^{-1}(\theta)\right)_{rr}\left(\SDM_x^{-1}(\theta)\right)_{r'r'}}
%d\theta \right)^{\frac{1}{2}}   \nonumber \\[3mm] &
%\stackrel{\eqref{equ_exp_inv_SDM_ensemble_i_i_d}}{=} 2\mincoh/(1+4\mincoh)
%\nonumber \\[3mm] & \stackrel{\mincoh\leq1/4}{\geq} \mincoh.\vspace{-2mm}
%\end{align}
\vspace{-2mm}
\begin{equation}
\rho_{x} \geq \mincoh \,. \vspace{-2mm}
\end{equation}
Therefore, the process $\mathbf{x}^{(i)}[\timeidx]$ belongs to $\processclass$
for any $\mincoh\leq 1/4$ and $\dmax \geq 1$.
The covariance matrix $\mathbf{C}_{i} = \expect \{ \mathrm{vec}(\xobs) \mathrm{vec}(\xobs)^{T} | i \}$ of
the observation $\xobs$, given the index $i$ satisfies
\vspace{-2mm}
\begin{equation}
\label{equ_expression_C_i_i_i_d}
\mathbf{C}_{i} = \mathbf{I}_{\samplesize} \otimes \mathbf{S}_{i}(0) \vspace{-2mm}
\end{equation} 
and, in turn,
\vspace{-2mm}
\begin{equation} 
\label{equ_expression_bar_C_i_i_d}
\overline{\mathbf{C}} = (1/M) \sum_{i \in[M]} \mathbf{C}_{i} =
\mathbf{I}_{\samplesize} \otimes \overline{\mathbf{S}} \vspace{-2mm}
\end{equation} 
with $\overline{\mathbf{S}} \defeq (1/M) \sum_{i \in [M]} \mathbf{S}_{i}(0)$.  
Inserting \eqref{equ_expression_C_i_i_i_d} and
\eqref{equ_expression_bar_C_i_i_d} into the bound
\eqref{equ_entropy_based_bound_MI},
%\begin{align}\label{32}
%I(\xobs;i) & \leq \log \big| \mathbf{I} \otimes \overline{\mathbf{S}}  \big| -
%(1/M) \sum_{i\in [M]} \log \big| \mathbf{I} \otimes \mathbf{S}_{i}(0) \big| \nonumber
%\\ & \stackrel{(a)}{=}  \samplesize \left( \log \big| \overline{\mathbf{S}} 
%\big| - (1/M) \sum_{i\in [M]} \log \big| \mathbf{S}_{i}(0) \big| \right),
%\end{align}
\vspace{-2mm}
\begin{equation}\label{32}
I(\xobs;i) \leq \samplesize \Big( \log \big| \overline{\mathbf{S}} 
\big| - (1/M) \sum_{i\in [M]} \log \big| \mathbf{S}_{i}(0) \big| \Big), \vspace{-2mm}
\end{equation}
where we used the identity $|\mathbf{I}_{\samplesize}
\otimes \mathbf{B}| =  |\mathbf{B}|^{\samplesize}$ (cf. \cite[Ch. 4]{Horn91}).
Note that the matrix $\mathbf{S}_{i}(0)$ has one eigenvalue equal to
$2(1-2a/(1+2a))$ and $p-1$ eigenvalues equal to $2$. Furthermore, setting $\gamma \defeq
1-\frac{2a}{(1+2a)p} + \frac{2a}{(1+2a)p(p-1)}$, the matrix
$\overline{\mathbf{S}}$ has one eigenvalue equal to
$2(\gamma-\frac{2a}{(p-1)(1+2a)})$ and $p-1$ eigenvalues equal to $2\gamma$.
Using the inequality $\log(1+x) \leq x$ for $x \geq 0$ and closely following the
calculation in \cite[Sec. 4.5.1.]{WeiWangPhd}, one obtains from \eqref{32}
\vspace{-2mm}	
\begin{equation} 
\label{equ_upper_bound_MI_i_i_d}
I(\xobs;i) \leq N 16 \mincoh^2. \vspace{-2mm} 
\end{equation}
Inserting \eqref{equ_upper_bound_MI_i_i_d} and
\eqref{equ_entropy_graph_ensemble_i_i_d} into
\eqref{equ_cond_MI_asmpt_reliable_sel} yields the bound
\eqref{necessary_cond_thm_iid}.
\vspace{-3mm}
\subsection{Proof of Theorem \ref{thm_processes}}\label{proof_processes}
\vspace{-1mm}
Consider the selection scheme $\widehat{\mathcal{N}}(r)$ described through
\eqref{decoder1}-\eqref{decoder3} for the neighborhood of a specific node $r$.
We aim at bounding the probability of failing to recover the correct neighborhood $\mathcal{N}(r)$,
$p_{err,r}=\Pr\left\{ \widehat{\mathcal{N}}(r) \neq \mathcal{N}(r) \right\}$.
We consider separately the two cases:
$\mathcal{N}(r)=\emptyset$ and $\mathcal{N}(r)\neq\emptyset$. \\
\noindent \textbf{ I. Node $r$ has no neighbor}\\
In this case, a selection error of the exhaustive search decoder can only occur
if for some node $r' \in V\setminus\left\{ r \right\}$ the statistic exceeds the
threshold, i.e., $|Z(r')| \geq \eta$.
Thus, $p_{a} = \Pr \left\{ \exists r' \in V\setminus\left\{ r \right\}
: |Z(r')| \geq \eta \right\} $ which, via a union bound, can be further bounded
as
\vspace{-2mm}
\begin{equation}
p_{a} \leq (p-1) \max_{r'\in V\setminus\left\{ r \right\}}
\Pr\left\{ |Z(r')| \geq \eta \right\}.\vspace{-2mm}
\end{equation}
\noindent
\textbf{II. Node $r$ has a single neighbor}\\
Suppose $(r,c)\in E$. The probability of erroneous detection of $\mathcal{N}(r)$
can be written as
\vspace{-2mm}
\begin{equation}
p_{b} = \Pr \left\{ \exists r' \in V \setminus
\left\{r,c\right\} : |Z(r')| \geq \eta \,\, \mathrm{and} \,\, |Z(c)| \leq \eta
\right\} \,,\vspace{-2mm}
\end{equation}
for some positive $\eta$. Using a union bound argument, this probability can be
bounded as
\vspace{-2mm}
\begin{equation}
p_{b} \leq (p-2) \Pr \left\{ \lvert Z(r') \rvert \geq \eta \right\} +
\Pr \left\{ \lvert Z(c) \rvert \leq \eta \right\} \,. \vspace{-2mm}
\end{equation}

Since $Z(\cdot)$ is the inner product of two zero-mean Gaussian vectors (cf.
\eqref{decoder1}), we can make direct use of \cite[Lemma E.2]{AJtechreport}
with $\mathbf{Q}=\mathbf{I}$ so separately bound $p_{a}$ and $p_{b}$ and in turn
$p_{err,r}$ since $p_{err,r} \leq \max\{p_{a},p_{b}\}$. In order to apply it to
above derivations, we need to characterize $\expect\left\{ Z(c) \right\} $ and
$\expect\left\{ Z(r') \right\} $ for $r'\notin\mathcal{N}(r)$.
Elementary calculations reveal that
\vspace{-2mm}
 \begin{equation}
 \expect\left\{ \hat{\mathbf{x}}_{r}^T \hat{\mathbf{x}}_{c}
\right\} = \left( \mathbf{C}_0 \right)_{rc} \sum_{\freqidx=1}^{2\samplesize} 
\left\lvert H\left(\frac{\freqidx}{2\samplesize}\right) \right\rvert ^2
=2\samplesize \left(
\mathbf{C}_0 \right)_{rc}
\vspace{-2mm}
 \end{equation}
 and $\frac{\mincoh}{B} \leq \left| \left( \mathbf{C}_0 \right)_{rc} \right|$.
 Further, $\expect\left\{ \hat{\mathbf{x}}_{r}^T \hat{\mathbf{x}}_{r'}
\right\} =0 $ for any $r' \notin \mathcal{N}(r)$. Based on \eqref{eq:AB}, we
can bound the spectral norm of the covariance matrix $\mathbf{C}_{\hat{x}_r} = E\left\{
 \hat{\mathbf{x}}_{r} \hat{\mathbf{x}}_{r}^{T}\right\}$ as
 $\|\mathbf{C}_{\hat{\mathbf{x}}_r}\|_2 \leq B \quad \forall r\in V$.
%  \vspace{-2mm}
%  \begin{equation}
%  \|\mathbf{C}_{\hat{\mathbf{x}}_r}\|_2 \leq B \quad \forall r\in V\,.
%  \vspace{-2mm}
%  \end{equation}
 Choosing the threshold in \eqref{decoder3} as $\nu = \mincoh/(2B)$, we obtain
 via \cite[Lemma E.2]{AJtechreport} the bound
 \vspace{-2mm}
\begin{equation}
p_{err,r} \leq 2p \exp \left(\frac{-N \mincoh^2}{32 B^4}\right). \vspace{-2mm}
\end{equation}
Using again a union bound argument, we obtain that the probability $p_{err}$ of
inferring an incorrect graph, that is, the probability of selecting at least one out of $p$ neighbor
sets incorrectly, is bounded as
\vspace{-2mm}
\begin{equation}
p_{err} \leq p\cdot p_{err,r} \leq 2p^2 \exp \left(\frac{-N \mincoh^2}{32
B^4}\right).\vspace{-2mm}
\end{equation}
Thus, if we require that the probability of a selection error does not exceed a
small number $\delta$, we obtain the sufficient condition \eqref{bound_processes}
on the sample-size $\samplesize$.

\vspace{-2mm}
\section{Conclusions}\label{sec:conclusions}
\vspace{-2mm}
We characterized the information theoretic limits of graphical model selection
for Gaussian time series by deriving a necessary condition on the sample-size
such that reliable selection may be possible. For a specific subclass of time
series with extremely sparse CIGs we showed that the necessary condition is
sharp.
In particular, we verified that a simple selection scheme is successful for a
sample-size close to the information theoretic limit. Somewhat unexpected, our
analysis reveals that in general the required sample-size is independent of the
correlation width, i.e., it does not depend on the amount of smoothness of the
SDM.
This suggests, in turn, that the sufficient condition presented in
\cite{JuHeck2014} for a novel selection scheme is far from optimal.

\renewcommand{\baselinestretch}{0.9}\normalsize\footnotesize
% \footnotesize
\bibliographystyle{IEEEbib}
%\bibliography{../../bibliography/tf-zentral,../../bibliography/work_valentin,../../bibliography/work_gabor}

\bibliography{work_gabor,tf-zentral,LitAJ_ITC}

\begin{thebibliography}{10}

\bibitem{Brockwell91}
P.~J. Brockwell and R.~A. Davis,
\newblock {\em Time Series: Theory and Methods},
\newblock Springer, New York, NY, 1991.

\bibitem{Luetkepol2005}
Helmut L{\"u}tkepohl,
\newblock {\em New Introduction to Multiple Time Series Analysis},
\newblock Springer, New York, 2005.

\bibitem{Dahlhaus2000}
R.~Dahlhaus,
\newblock ``Graphical interaction models for multivariate time series,''
\newblock {\em Metrika}, vol. 51, pp. 151--172, 2000.

\bibitem{WangWain2010}
W.~Wang, M.J. Wainwright, and K.~Ramchandran,
\newblock ``Information-theoretic bounds on model selection for gaussian markov
  random field,''
\newblock in {\em Proc. IEEE ISIT-2010}, Austin, TX, Jun. 2010, pp. 1373--1377.

\bibitem{PetersCausalInf2012}
J.~{Peters}, D.~{Janzing}, and B.~{Sch{\"o}lkopf},
\newblock ``Causal inference on time series using structural equation models,''
\newblock {\em ArXiv e-prints}, Jul. 2012.

\bibitem{Bento2010}
J.~Bento, M.~Ibrahimi, and A.~Montanari,
\newblock ``Learning networks of stochastic differential equations,''
\newblock in {\em Advances in Neural Information Processing Systems 23},
  Vancouver, CN, 2010, pp. 172--180.

\bibitem{Songsiri09}
J.~Songsiri, J.~Dahl, and L.~Vandenberghe,
\newblock ``Graphical models of autoregressive processes,''
\newblock in {\em Convex Optimization in Signal Processing and Communications},
  Y.~Eldar and D.~Palomar, Eds., pp. 89--116. Cambridge Univ. Press, Cambridge,
  UK, 2010.

\bibitem{songsiri2010}
J.~Songsiri and L.~Vandenberghe,
\newblock ``Topology selection in graphical models of autoregressive
  processes,''
\newblock {\em Journal of Machine Learning Research}, vol. 11, pp. 2671--2705,
  2010.

\bibitem{JuHeck2014}
A.~Jung, R.~Heckel, H.~B\"{o}lcskei, and F.~Hlawatsch,
\newblock ``Compressive nonparametric graphical model selection for stationary
  time series,''
\newblock in {\em Proc. IEEE ICASSP 2014}, Florence, Italy, May 2014.

\bibitem{LauritzenGM}
S.~L. Lauritzen,
\newblock {\em Graphical Models},
\newblock Clarendon Press, Oxford, UK, 1996.

\bibitem{Brillinger96remarksconcerning}
R.~Brillinger,
\newblock ``Remarks concerning graphical models for time series and point
  processes,''
\newblock {\em Revista de Econometria}, vol. 16, pp. 1--23, 1996.

\bibitem{Wain2009TIT}
M.~J. Wainwright,
\newblock ``Information-theoretic {L}imits on {S}parsity {R}ecovery in the
  {H}igh-{D}imensional and {N}oisy {S}etting,''
\newblock {\em IEEE Trans. Inf. Theory}, vol. 55, no. 12, pp. 5728--5741, Jun.
  2009.

\bibitem{CaiZhaoZhou2013}
T.~T. Cai, Z.~Ren, and H.~H. Zhou,
\newblock ``Optimal rates of convergence for estimating toeplitz covariance
  matrices,''
\newblock {\em Probability Theory and Related Fields}, vol. 156, no. 1-2, pp.
  101--143, 2013.

\bibitem{InverseBCM83}
T.~De Mazancourt and D.~Gerlic,
\newblock ``The inverse of a block-circulant matrix,''
\newblock {\em IEEE Trans. Antennas and Propagation}, vol. AP-31, no. 5, Sep.
  1983.

\bibitem{coverthomas}
T.~M. Cover and J.~A. Thomas,
\newblock {\em Elements of Information Theory},
\newblock Wiley, New Jersey, 2 edition, 2006.

\bibitem{WeiWangPhd}
W.~Wang,
\newblock {\em Sparse signal recovery using sparse random projections},
\newblock Ph.D. thesis, EECS Department, University of California, Berkeley,
  Dec 2009.

\bibitem{HendersonSearle}
H.~V. Henderson and S.~R. Searle,
\newblock ``On deriving the inverse of a sum of matrices,''
\newblock {\em SIAM Review}, vol. 23, no. 1, pp. 53--60, Jan. 1981.

\bibitem{Horn91}
Roger~A. Horn and Charles~R. Johnson,
\newblock {\em Topics in {M}atrix {A}nalysis},
\newblock Cambridge Univ. Press, Cambridge, UK, 1991.

\bibitem{AJtechreport}
A.~Jung, R.~Heckel, H.~B\"{o}lcskei, and F.~Hlawatsch,
\newblock ``Compressive nonparametric graphical model selection for stationary
  time series: A multitask learning approach,''
\newblock Tech. {R}ep., Vienna University of Technology, Institute of
  Telefommunications, March 2014.

\end{thebibliography}

\end{document}